\documentclass[11pt]{article}
\usepackage{xspace,amssymb,amsmath,amscd,amsthm,epsfig}

\advance \topmargin by -\headheight
\advance \topmargin by -\headsep
     
\evensidemargin \oddsidemargin
\theoremstyle{plain}
\newtheorem{theorem}{Theorem}

\newtheorem{conjecture}[theorem]{Conjecture}
\theoremstyle{definition}

\newtheorem{example}[theorem]{Example}
\newtheorem{remark}[theorem]{Remark}
\newtheorem{problem}[theorem]{Problem}

\title{A combinatorial version of Sylvester's four-point problem}
\author{Gregory S. Warrington\\
  gwarring@cems.uvm.edu}

\begin{document}

\maketitle

\begin{abstract}
J. J. Sylvester's four-point problem asks for the probability that
four points chosen uniformly at random in the plane have a triangle as
their convex hull.  Using a combinatorial classification of points in
the plane due to Goodman and Pollack, we generalize Sylvester's
problem to one involving reduced expressions for the long word in
$S_n$.  We conjecture an answer of $1/4$ for this new version of the
problem.
\end{abstract}

\begin{remark}
  O. Angel and A. E. Holroyd~\cite{holroyd} prove a more general
  result which implies our Conjecture~\ref{mainconj}.
\end{remark}
\section{Introduction}

Fix $n \geq 4$.  Pick uniformly at random a reduced expression
$\mathsf{w}$ for the long word in $S_n$.  Then pick uniformly at
random a $4$-subset of $\{1,\ldots,n\}$.  This pair of choices
determines a reduced expression $\mathsf{v}$ for the long word in
$S_4$.  What is the probability that
\begin{equation}\label{eq:reent}
  \mathsf{v} \in X := \{s_1s_2s_3s_2s_1s_2,\ s_3s_2s_1s_2s_3s_2,\
  s_2s_1s_2s_3s_2s_1,\ s_2s_3s_2s_1s_2s_3\}?
\end{equation}

\begin{conjecture}\label{mainconj}
  For any $n\geq 4$, the probability is $1/4$.
\end{conjecture}

For $n=4$, the probability of $1/4$ can be computed directly from the
possible $\mathsf{w}$ listed in Figure~\ref{fig:all16}.  The
conjecture can be checked for small values of $n$ by computer:
\begin{eqnarray*}
  {\text{5 points}} & \frac{960}{\binom{5}{4} 768} &= \frac{1}{4}\\
  {\text{6 points}} & \frac{1,098,240}{\binom{6}{4} 292,864} &= \frac{1}{4}\\
  {\text{7 points}} & \frac{9,631,498,240}{\binom{7}{4} 1,100,742,656} &= \frac{1}{4}\\
  {\text{8 points}} & \frac{850,653,924,556,800}{\binom{8}{4} 48,608,795,688,960} &= \frac{1}{4}.
\end{eqnarray*}
The computer check for $n=8$ took several thousand hours on 3GHz CPUs.

\begin{remark}
  The behavior of arbitrary reduced expressions for the long word in
  $S_n$ (also known as \emph{sorting networks}) has been considered in
  a number of contexts.  Most notably, Angel et al.~\cite{AHRV}
  consider several convergence questions.  V. Reiner~\cite{vic}
  computes the expected number of possible Yang-Baxter moves for such
  a reduced expression in the symmetric group while
  B. Tenner~\cite{bridget} performs the analogous calculation for the
  hyperoctahedral group.  The above results all use connections to the
  theory of Young tableaux (see~\cite{E-G,rstan}).  Through the
  Goodman-Pollack correspondence~\cite{G-P}, described here in
  Section~\ref{sec:mot}, there are also connections to halving lines
  and $k$-sets of points in the plane (see, e.g.,~\cite{lovasz,toth}).
  The nature of the connections between Conjecture~\ref{mainconj} and
  these other works is still unclear.
\end{remark}

\section{Motivation}
\label{sec:mot}

Conjecture~\ref{mainconj} combines a question of J.J. Sylvester with a
combinatorial classification of points in the plane due to Goodman and
Pollack.

In 1864, J.J. Sylvester~\cite{Sylv} posed the following
\begin{problem}\label{prob:orig}
  Given four points chosen at random in the plane, find the probability
  that they form a reentrant (rather than convex) quadrilateral.
\end{problem}
As was quickly realized at the time, the problem as stated is
ill-posed; it does not specify how the plane is being modeled.
Woolhouse (see~\cite{Woolh}) realized that it sufficed to pick points
in a closed, bounded region.  The probability would be left invariant
under the scaling necessary to model the plane.  In fact, the exact
probability is computable by integration for any convex region.  Among
convex regions, a triangle maximizes the probability at a value of
$1/3$ and a disk minimizes it at $35\pi^2/12\approx 0.295$.  The
reader is referred to Pfiefer~\cite{Pfiefer} for an excellent history
of the problem.

Less than twenty years later, R. Perrin~\cite{Perrin} considered the
matter of combinatorially classifying collections of $n$ points in
general position in the plane.  More recently, this work has been
extended by Goodman and Pollack~\cite{G-P} in the form of
\emph{allowable sequences} (see also Knuth~\cite[\S 8]{knuth}).  There
are several variations on Goodman and Pollack's allowable sequences
one can work with.  Our setup is as follows.

For $n \geq 1$, let $\mathcal{P}_n$ denote the set of all possible
$n$-tuples of points in the plane in general position.  We consider
the symmetric group $S_n$ on $\{1,2,\ldots,n\}$ as generated by the
adjacent transpositions $s_i = (i,i+1)$ for $1\leq i < n$.  For any
$w\in S_n$, a \emph{reduced expression} for $w$ is a word
in the $s_i$'s of minimum length among all those whose product equals
$w$.  Given $w\in S_n$, let $\mathcal{R}(w)$ to be the set of all
possible reduced expressions for $w$.  For the long word $w_0 =
[n,n-1,\ldots,2,1]$, we now define a map $\phi: \mathcal{P}_n
\rightarrow \mathcal{R}(w_0)$.  The map $\phi$ depends on a fixed
directed line $\ell$ in the plane which can be chosen arbitrarily.
Denote by $\ell(\theta)$ the line $\ell$ rotated counterclockwise
through $\theta$ radians about some point on $\ell$.

Pick $P \in \mathcal{P}_n$.  For each pair of distinct points $i,j\in
P$, there is a unique angle $\theta' \in [0,\pi)$ for which
$\ell(\theta')$ is orthogonal to the line through $i$ and $j$.
Running over all unordered pairs of points, these $\theta'$ define a
sequence
\begin{equation}\label{eq:thetaseq}
  0\leq \theta_1 < \theta_2 < \cdots <
  \theta_{\binom{n}{2}} < \pi.
\end{equation}
We now label the points according to their projections onto $\ell$.
(In the probability-zero case that $\theta_1 = 0$, we instead project
onto $\ell(-\epsilon)$ for some sufficiently small $\epsilon$.)
For each angle $\theta$ not occurring in the sequence
\eqref{eq:thetaseq}, the projection of $P$ onto $\ell(\theta)$
determines a permutation $\pi_\theta$ with respect to the initial
labeling.  Write $\pi_{\theta_0} = [1,2,\ldots,n]$.  We set $\phi(P)
:= \mathsf{w} = \mathsf{w_1w_2\cdots w_{\binom{n}{2}}}$ where $w_k =
s_i$ if $\pi_{\theta_k+\delta} = \pi_{\theta_{k-1}}s_i$ for $\delta =
(\theta_{k+1}-\theta_k)/2$.

\begin{figure}[htbp]
    \centering
      {\scalebox{.55}{\includegraphics{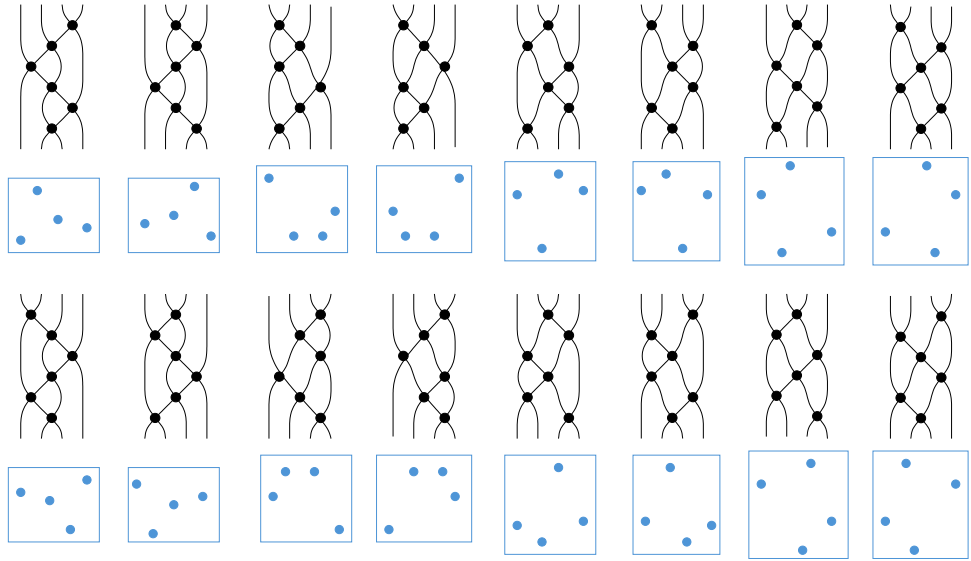}}}
      \caption{Reduced expressions for $w_0\in S_4$.}
      \label{fig:all16}
\end{figure}

We illustrate $\phi(\mathcal{P}_4)$ in Figure~\ref{fig:all16}.  Each
expression $\mathsf{w} \in \mathcal{R}(w_0)$ is written as a string
diagram (see, e.g., \cite{321-av}).  To clarify conventions, we note
that the first reduced expression in the first row of the figure is
$s_2s_3s_2s_1s_2s_3$.  Below each $\mathsf{w}$ is a configuration of
points $P$ for which $\phi(P) = \mathsf{w}$ (assuming $\ell$ is
horizontal and directed to the right).  It is worth pointing out that
Goodman and Pollack~\cite{G-P} show that $\phi$ is not surjective for
$n\geq 5$.

The reentrant configurations in the first two columns of
Figure~\ref{fig:all16} correspond to the reduced expressions
comprising the set $X$ in~\eqref{eq:reent}.

\section{Combinatorial Version}

\begin{figure}[htbp]
    \centering
      {\scalebox{.5}{\includegraphics{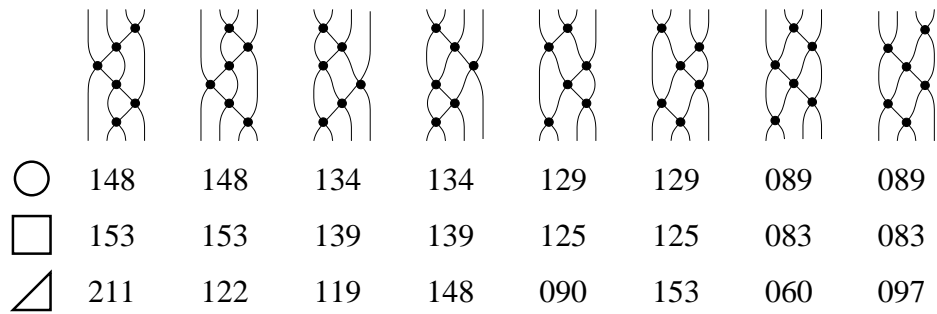}}}
      \caption{Monte Carlo simulations for three regions (1000 trials).}
      \label{fig:monte}
\end{figure}

We now rephrase Sylvester's problem in terms of the elements of
$\mathcal{R}(w_0)$.  Fix a closed, bounded region $R$ in the plane
along with a directed line $\ell$.  The uniform distribution on $R$
induces a probability distribution $f_n$ on $\phi(\mathcal{P}_n)$ for
each $n$.  We can think of picking four points at random from $R$ as
picking one of the 16 elements of $\mathcal{R}([4,3,2,1])$ according
to $f_n$.  The distribution $f_n$ is displayed for several different
shapes in Figure~\ref{fig:monte}.  For brevity in the figure, we have
identified the two rows of Figure~\ref{fig:all16}.  (Entries in
the same column of Figure~\ref{fig:all16} correspond to choosing
the opposite direction for the line $\ell$ and will therefore appear
with equal frequencies in $f_n$.)

For the general, combinatorial version of this problem, we dispense
with the geometry of our region $R$.  We start with the \emph{uniform}
distribution on $\mathcal{R}(w_0)$ rather than the geometrically
induced distribution $f_n$.

\begin{problem}\label{prob:comb}
  For $n\geq 4$ and $w_0\in S_n$, pick $\mathsf{w}\in\mathcal{R}(w_0)$
  uniformly at random and a four-subset $\{a,b,c,d\} \subseteq
  \{1,2,\ldots,n\}$ uniformly at random.  The pair of choices induces
  a reduced expression $\mathsf{v}$ for the long word in
  $S_{\{a,b,c,d\}} \cong S_4$.  Find the probability that $\mathsf{v}\in X$.
\end{problem}

\begin{figure}[htbp]
    \centering
      {\scalebox{.35}{\includegraphics{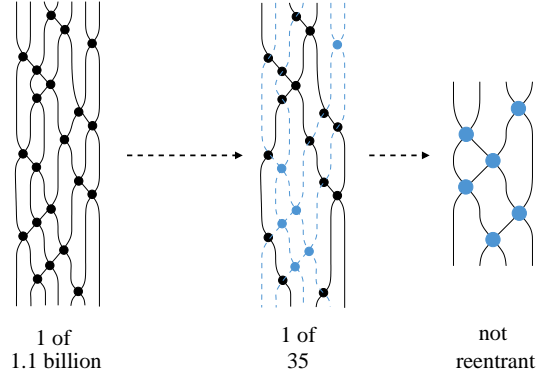}}}
      \caption{Two-stage combinatorial version illustrated for $n=7$.}
      \label{fig:7ptex}
\end{figure}

An example of the procedure of Problem~\ref{prob:comb} is illustrated
in Figure~\ref{fig:7ptex}.

\begin{example}
  R. Stanley~\cite{rstan} showed that $|\mathcal{R}(w_0)|$ equals the
  number of standard Young tableaux of staircase shape
  $(n-1,n-2,\ldots,1)$.  Applying the Frame-Robinson-Thrall hook
  formula for the number of such tableaux, we find that there are
  $768$ such reduced expressions when $n=5$.  They fall into three
  classes according to how many of the $\binom{5}{4}$ four-tuples lead
  to $\mathsf{v}$ being in $X$.  Figure~\ref{fig:5pts-classes}
  displays the size of each class, a representative, and a point
  configuration realizing one of the reduced expressions in the class.
  (Recall that some of the expressions are not realizable; the point
  configurations are for intuition only.)
  
\begin{figure}[htbp]
  \centering
  {\scalebox{.4}{\includegraphics{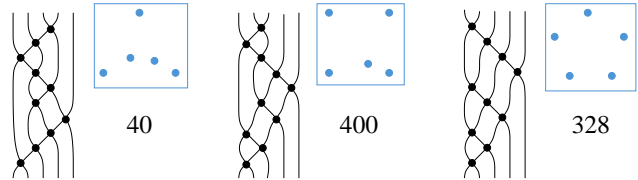}}}
  \caption{Computation of probability for $n=5$.}\label{fig:5pts-classes}
\end{figure}
From this data, we can calculate
\begin{equation}
  \mathrm{Pr}(\mathsf{v}\ \text{reentrant}) = 
  \frac{40\cdot 4 + 400\cdot 2 + 328\cdot 0}{\binom{5}{4}(40 + 400 + 328)} = \frac{1}{4}.
\end{equation}
Similar computations can be made on computer for $n \in \{6,7,8\}$.

\end{example}

\begin{remark}
  As can be seen in Figure~\ref{fig:all16}, the reduced expressions
  fall into four natural classes.  In terms of points in the plane,
  these are generated by the direction of rotation (clockwise or
  counterclockwise) and the choice of direction for the directed line
  $\ell$.  A natural assumption is that Conjecture~\ref{mainconj}
  holds for the three ``convex'' classes as well.  This is not true.
  For example, with $n=5$, the number of pairs
  $(\mathsf{w},\mathsf{v})$ falling into the each of the four classes
  is, respectively, $960$, $980$, $972$ and $928$.
\end{remark}



\bibliography{mybib}{}

\begin{thebibliography}{10}

\bibitem{holroyd}
Omer Angel and Alexander~E. Holroyd.
\newblock Random subnetworks of random sorting networks.
\newblock {\em Electron. J. Combin.}, 17(1):Note 23, 7, 2010.

\bibitem{AHRV}
Omer Angel, Alexander~E. Holroyd, Dan Romik, and B{\'a}lint Vir{\'a}g.
\newblock Random sorting networks.
\newblock {\em Adv. Math.}, 215(2):839--868, 2007.

\bibitem{321-av}
Sara~C. Billey and Gregory~S. Warrington.
\newblock Kazhdan-{L}usztig polynomials for 321-hexagon-avoiding permutations.
\newblock {\em J. Algebraic Combin.}, 13(2):111--136, 2001.

\bibitem{E-G}
Paul Edelman and Curtis Greene.
\newblock Balanced tableaux.
\newblock {\em Adv. in Math.}, 63(1):42--99, 1987.

\bibitem{G-P}
Jacob~E. Goodman and Richard Pollack.
\newblock On the combinatorial classification of nondegenerate configurations
  in the plane.
\newblock {\em J. Combin. Theory Ser. A}, 29(2):220--235, 1980.

\bibitem{Woolh}
C.~M. Ingleby.
\newblock In W.~J.~C. Miller, editor, {\em Mathematical questions and their
  solutions from the Educational Times}, volume~5, page 108. April 1865.

\bibitem{knuth}
D.~E. Knuth.
\newblock {\em Axioms and hulls}, volume 606 of {\em Lecture Notes in Computer
  Science}.
\newblock Springer-Verlag, Berlin, 1992.

\bibitem{lovasz}
L.~Lov{\'a}sz.
\newblock On the number of halving lines.
\newblock {\em Ann. Univ. Sci. Budapest. E\"otv\"os Sect. Math.}, 14:107--108,
  1971.

\bibitem{Perrin}
Perrin.
\newblock Sur le probl\`eme des aspects.
\newblock {\em Bull. Soc. Math. France}, 10:103--127, 1882.

\bibitem{Pfiefer}
Richard~E. Pfiefer.
\newblock The historical development of {J}. {J}. {S}ylvester's four point
  problem.
\newblock {\em Math. Mag.}, 62(5):309--317, 1989.

\bibitem{vic}
Victor Reiner.
\newblock Note on the expected number of {Y}ang-{B}axter moves applicable to
  reduced decompositions.
\newblock {\em European J. Combin.}, 26(6):1019--1021, 2005.

\bibitem{rstan}
Richard~P. Stanley.
\newblock On the number of reduced decompositions of elements of {C}oxeter
  groups.
\newblock {\em European J. Combin.}, 5(4):359--372, 1984.

\bibitem{Sylv}
J.~J. Sylvester.
\newblock {\em The Educational Times}, April, 1864.

\bibitem{bridget}
Bridget~Eileen Tenner.
\newblock On expected factors in reduced decompositions in type {$B$}.
\newblock {\em European J. Combin.}, 28(4):1144--1151, 2007.

\bibitem{toth}
G.~T{\'o}th.
\newblock Point sets with many {$k$}-sets.
\newblock {\em Discrete Comput. Geom.}, 26(2):187--194, 2001.
\newblock {ACM} Symposium on Computational Geometry (Hong Kong, 2000).

\end{thebibliography}
\bibliographystyle{plain}

\end{document}